\documentclass[journal]{IEEEtran}
\usepackage{amsmath}
\usepackage{}
\usepackage{amssymb}
\usepackage{graphicx}
\usepackage{colortbl}

\usepackage{algorithm}
\usepackage{algorithmic}

\usepackage{amsfonts}
\usepackage{latexsym,amsmath,amssymb}
\usepackage{graphicx,psfrag}

\newcommand{\no}{\nonumber}
\newcommand{\be}{\begin{equation}}
\newcommand{\ee}{\end{equation}}
\def\bee#1\eee{\begin{align}#1\end{align}}
\def\beee#1\eeee{\begin{align*}#1\end{align*}}
\def\le{\left}
\def\ri{\right}
\def\no{\nonumber}
\def\bee#1\eee{\begin{align}#1\end{align}}
\def\beee#1\eeee{\begin{align*}#1\end{align*}}
\newcommand{\bea}{\begin{eqnarray}}
\newcommand{\eea}{\end{eqnarray}}
\def\ba{\begin{array}}
\def\ea{\end{array}}

\def\ok{$\blacksquare$}

\newcommand{\bse}{\begin{subequations}}
\newcommand{\ese}{\end{subequations}}
\newtheorem{theorem}{Theorem}%[section]
\newtheorem{lemma}{Lemma}%[section]
\newtheorem{remark}{Remark}%[section]
\newtheorem{definition}{Definition}%[section]
\newtheorem{corollary}{Corollary}%[section]

\newtheorem{example}{Example}
\ifCLASSINFOpdf
\else
\fi
\begin{document}
\title{Criteria for stabilizing a multi-delay stochastic system with multiplicative control-dependent noises}

\author{
Cheng Tan \IEEEmembership{Member,~IEEE},
%Haoting Sui,
%Yuzhe Li,
Zhengqiang Zhang \IEEEmembership{Senior Member,~IEEE},
Haoting Sui,\\
Wing Shing Wong  \IEEEmembership{Life Fellow,~IEEE}
           % <-this % stops a space
\thanks{This work was supported in part by the National Natural Science Foundation of China under Grants 62173206;
the Natural Science Foundation of Shandong Province under Grant ZR2021ZD13;
China Postdoctoral Science Foundation under Grant 2021M691849.}

\thanks{C. Tan, Z. Zhang and H. Sui are with the School of Engineering, QuFu Normal University, Rizhao, Shandong 276800, China,
 (e-mail: tancheng1987love@163.com, qufuzzq@126.com).}

%\thanks{Y. Li is with the State Key Laboratory of Synthetical Automation for Process Industries, Northeastern University, Shenyang 110004, China.
%%the College of Electrical Engineering and Automation Shandong University of Science and Technology, Qingdao 266590, China (e-mail: hszhang@sdu.edu.cn).
%}

\thanks{
W. S. Wong is with the Department of Information Engineering, The Chinese University of Hong Kong, Shatin, N. T., Hong Kong (e-mail: wswong@ie.cuhk.edu.hk).}
}

%\thanks{
%H. Zhang is with the School of Control Science and Engineering, Shandong University, Jinan, Shandong 250061, China (\small  e-mail: hszhang@sdu.edu.cn). }

%\thanks{W. S. Wong is with the Department of Information Engineering, The Chinese University of Hong Kong, Shatin, N. T., Hong Kong (\small  e-mail: wswong@ie.cuhk.edu.hk).}


%\markboth{Journal of \LaTeX\ Class Files,~Vol.~11, No.~4, December~2012}%
%{Shell \MakeLowercase{\textit{et al.}}: Bare Demo of IEEEtran.cls for Journals}



\maketitle

\begin{abstract}
In this paper, we investigate the mean-square stabilization for discrete-time stochastic systems that endure both multiple input delays and multiplicative control-dependent noises.
For such multi-delay stochastic systems, we for the first time put forward two stabilization criteria: Riccati type and Lyapunov type.
On the one hand, we adopt a reduction method to reformulate the original multi-delay stochastic system to a delay-free auxiliary system and present their equivalent proposition for stabilization.
Then, by introducing a delay-dependent algebraic Riccati equation (DDARE), we prove that the system under consideration is stabilizable if and only if the developed DDARE has a unique positive definite solution.
On the other hand, we characterize the delay-dependent Lyapunov equation (DDLE)-based criterion, which can be verified by linear matrix inequality (LMI) feasibility test.
Besides, under some restricted structure, we propose an existence theorem of delay margin and more importantly, derive an explicit formula for computing its exact value.
\end{abstract}


% \fi
\IEEEpeerreviewmaketitle

\section{Introduction}

In the last decades, the stability/stabilization issues for stochastic systems have attracted considerable interest because of their extensive applications
in economics as well as in engineering fields; See \cite{A2012}--\cite{deng2022} and the references therein. %for a partial list of references.
For linear time invariant systems, a large body of excellent results, including necessary and sufficient conditions, have been developed.
In particular, the Lyapunov-type stabilization criteria was derived in terms of the feasibility of a certain LMI in \cite{zhouxy2000}, while the Riccati-type result was developed via a unique positive definite solution satisfying generalized algebraic Riccati equation (GARE) in \cite{huang2008}.
As a supplement, the mean-square stabilization was first characterized by the spectral locations of coefficient matrices in \cite{zhangwh2004}.

The aforementioned studies were exclusively concerned with delay-free stochastic models.
Recently, as the product of the integration of control technology and network communication technology, networked control systems (NCSs) have been actively researched due to their widespread applications.
Because of the limitation of network bandwidth power, the signal transmitted in unreliable channel may encounter various network-induced uncertainties including data packet loss (multiplicative noise) and transmission delay.
As a result, some studies have concentrated on the scenarios where time delay and random noises occur simultaneously \cite{yued2005}-\cite{tan2019}.
The most popular methods, such as Lyapunov-Krasivskii functional approach, are mainly suited for exploiting sufficient stabilization conditions by virtue of LMIs.
Specifically, for single-delay stochastic model, the necessary and sufficient stabilization condition was first presented in terms of coupled algebraic Riccati equations (CAREs) in \cite{lilin2015}.
However, it is difficult to calculate the positive definite solutions because CAREs are nonlinear and their positive solutions are coupled.
How to best utilize available information to design a stabilizing control policy as well as the search for stabilization criteria for more general multi-delay stochastic dynamics remain open and challenging questions.

On the other hand, delay margin, as a fundamental measure of robust stabilization against uncertain delays, has also received extensively attention; See \cite{margin01}--\cite{margin03}.
There are two natural questions for multi-delay stochastic systems: What is the largest range of delay so that there exists an admissible control law that can stabilize the target system within that entire range? How to calculate this delay margin?
Unfortunately, to the best of our knowledge, there does not exist work on characterizing the delay margin for multi-delay stochastic model, which motivates us to undertake an in-depth study.

Facilitated by stochastic control techniques, we seek to provide a series of stabilization criteria for stochastic system with multiple input delays and multiplicative noises.
Different from the previous work, one significant contribution is that our control law is designed as the feedback of an extended state that contains the recent available state information and part values from previous control inputs.
It is remarkable that the developed criteria are necessary and sufficient, which are first obtained within the framework of multi-delay stochastic system with multiplicative noises in control variables.
These criteria run in parallel to the classical results in \cite{zhouxy2000,huang2008}.

Our research methodology is described as follows.
First, motivated by predictive control methodology proposed in \cite{art1982}, we adopt a reduction strategy to transform the original multi-delay stochastic system into an auxiliary delay-free model and demonstrate their equivalent proposition for stabilization.
Then, we present the Riccati-type stabilization criterion and the design procedure of the stabilizing control law.
The expression is that system is stabilizable if and only if the predefined DDARE has a unique positive definite solution.
Of equal importance, we characterize the Lyapunov-type stabilization criterion by means of DDLE.
Utilizing Schur complement decomposition and matrix transformation technique, this criterion can be expressed by LMI-based condition.
Specifically, as an application, we apply our theory to study the delay margin problem.
On the basis of the stabilization criteria and operator theory, we show that there exists a unique delay margin for some restricted single-delay stochastic model.
More importantly, some explicit formula for computing the delay margin is derived for uncoupled system.
At last, we present two simple examples to confirm
%show the effectiveness of
our theoretic results.

Notation: For any integer $i< j$, we define $\mathbb{N}_{[i,j]}\triangleq \{i,\cdots,j\}$.
$Z\geq0~(>0)$ means that $Z$ is a positive semidefinite (positive definite) matrix, and $Z_1\geq Z_2~(>Z_2)$ means that $Z_1-Z_2\geq0~(>0)$.
$\{\omega_t,t\in\mathbb{N}\}$ denotes a sequence of real random variables defined on the complete filtered probability space $(\Omega,\mathcal{F},\mathcal{P};\mathcal{F}_t)$ with $\mathcal{F}_t=\sigma\{\omega_s,~s\in\mathbb{N}_{[0,t]} \}$.
Define $\hat{x}_{s|t}={\bf E}(x_s|\mathcal{F}_t)$ which signifies the conditional expectation of $x_s$ w.r.t. $\mathcal{F}_t$.

\section{Problem formulation}
In this paper, we consider the following discrete-time stochastic system with both multiple delays and multiplicative control-dependent noises
\be
x_{t+1}=Ax_t+\sum_{\tau=0}^D \le(B_\tau +\omega^\tau_t C_\tau \ri)u_{t-\tau},\label{system01}
\ee
where $x_t\in\mathbb{R}^n$ is the state, and $u_t\in\mathbb{R}^m$ is the control input executed at time $t\geq 0$.
System matrices
%$A$, $B_\tau$, and $C_\tau,~\tau\in\mathbb{N}_{[0,D]}$, are constant of compatible dimensions, where $\mathbb{N}_{[i,j]}\triangleq \{i,\cdots,j\}$.
Assume that the initial conditions $x_0$, $u_t$, $t\in \mathbb{N}_{[-d,-1]}$, are given a priori.
The multiplicative noises are assumed to be random sequences with independent and identically distributed realizations
\be
{\bf E}\le(\omega_t^\tau \ri)=0,
~{\bf E}\le(\omega_t^\tau \omega_s^\tau\ri)=\sigma_\tau^2\delta_{ts}, ~\forall\tau\in\mathbb{N}_{[0,D]}
\ee
where $\delta_{ts}$ is a Kronecker function.
%and $\mathbb{N}_{[i,j]}\triangleq \{i,\cdots,j\}$.
We further denote
\be
\omega_t\triangleq \left(\omega_t^0~\omega_t^1\cdots \omega_t^D \right)',~t\geq0.
\ee
which defines the $\sigma$-algebra as $\mathcal{F}_t=\sigma\{\omega_s,~s\in\mathbb{N}_{[0,t]} \}$.
%To facilitate the narrative, define $\hat{x}_{s|t}={\bf E}(x_s|\mathcal{F}_t)$ which signifies the conditional expectation of $x_s$ w.r.t. $\mathcal{F}_t$.
%For notational convenience, we describe system (\ref{system01}) as $\left(A,B,C,T\right)$.
Before proceeding further, we impose the following definition of asymptotical mean-square stabilization.
\begin{definition}
\label{def01}
System \eqref{system01} is said to be asymptotically mean-square stabilizable, if there exists a feedback control input $u_t$
such that the corresponding closed-loop system  is asymptotically mean-square stable, i.e. for any initial values $x_0$ and $u_t,~t\in\mathbb{N}_{[-D,-1]}$, the state $x_{t}$ in (\ref{system01}) satisfies $\lim_{t\rightarrow \infty} {\bf E}(x_{t}'x_t)=0$.
\end{definition}

This paper is concerned with the asymptotical mean-square stabilization problem for multi-delay stochastic system \eqref{system01} in which the control input should be $\mathcal{F}_{t-1}$-measurable.
It is remarkable that if the control law is designed to be a feedback of the state, owing to information gap, it is difficult to obtain the necessary and sufficient condition \cite{lilin2015,tan2019}.
%Until now, to the best of our knowledge, there have been no firm results to this problem.
To tackle this problem, motivated by Smith predictor in \cite{Ito1981}, we introduce a viable controller construction by utilizing an extended state vector that contains the recent available state and previous control inputs in finite horizon.
In this paper, the set of admissible control input is given as
\be
\mathcal {U}_{ad}\triangleq \le\{u_t\in L^2_{\mathcal{F}},u_{t}=K_0 x_t+\sum_{\tau=1}^D K_\tau u_{t-\tau} \ri\},\label{ad-input-set}
\ee
where
\bee
L^2_{\mathcal{F}}\triangleq \{ u(t)~is~\mathcal{F}_{t}-measurable,\sum_{t=0}^{\infty} {\bf E}(u_t'u_t)<\infty \}. \label{ad-input-set00}
\eee

The problems to be solved are formulated as follows.

\begin{itemize}
  \item Explore some control strategy $u_t\in\mathcal {U}_{ad}$ to stabilize system (\ref{system01}), while develop the necessary and sufficient stabilization conditions.
  \item Explore exact delay margin to guarantee stabilization.
\end{itemize}
%
%\noindent{\bf Problem 1} Explore some control strategy $u_t\in\mathcal {U}_{ad}$ to stabilize system (\ref{system01}), while develop the necessary and sufficient stabilization conditions.
%
%\noindent{\bf Problem 2} Explore exact delay margin to guarantee stabilization.

\begin{figure}[thpb]
      \centering
      \includegraphics[scale=0.33]{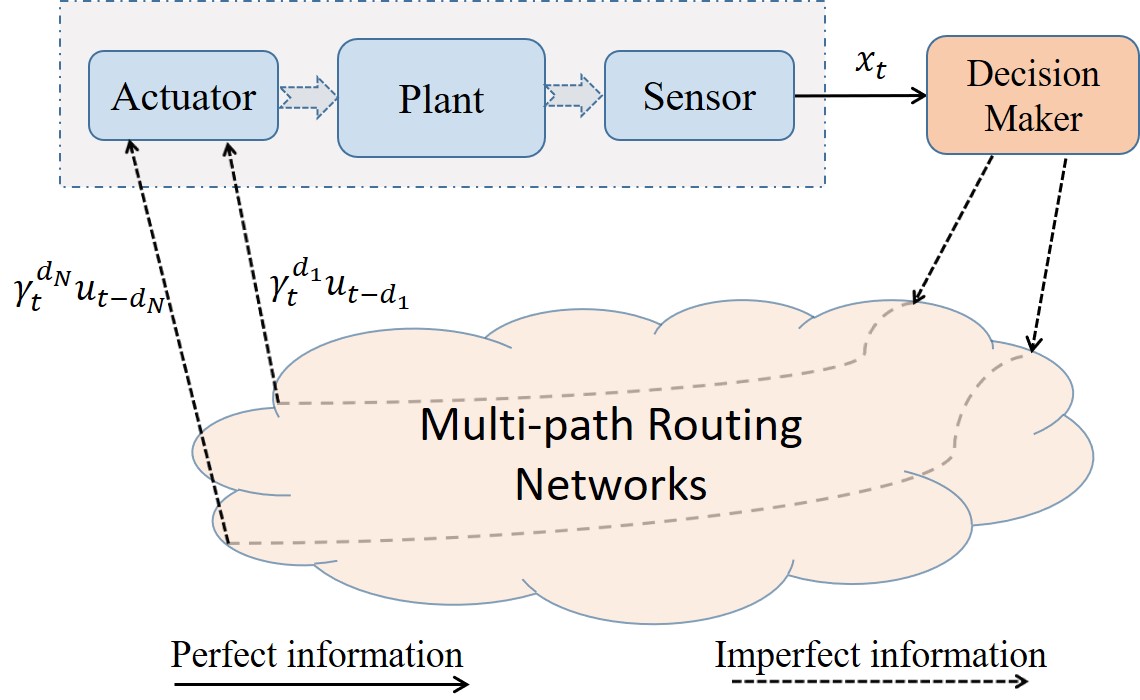}
      \centering\caption{\scriptsize{WNCS over multi-path routing network}}
\end{figure}

\begin{remark}
Note that the considered stochastic system has a wide application in wireless NCS (WNCS).
Notice that the considered WNCS with both transmission delay and packet loss in \cite{tan2015} is a special case of system \eqref{system01}.
More generally, as a flexible driving architecture in wireless sensor and Ad-hoc networks, the multi-path routing approach has gained popularity for various network management purposes \cite{multi01}--\cite{multi03}.
Fig. 1 shows a typical structure of WNCS over unreliable multi-path routing network, which is responsible for the transmission of control signals, where the delivered control input over the $i$-th path suffers both constant delay $d_i\geq 0$ and data packet loss.
Here, the arrival process of the control signal is modeled via a binary random variable $\gamma^{i}_t$ defined as
\be
\gamma^{i}_t=
\le\{
\begin{aligned}
& 1, ~\text{if $u_{t-d_i}$ has been delivered to the actuator},\\
& 0 , ~\text{otherwise}.
\end{aligned}
\ri.
\ee
Therefore, the dynamics of the overall WNCS follows
\be
x_{t+1}=Ax_t+\sum_{i=1}^N \gamma^{i}_t Bu_{t-d_i}.\label{NCS-original}
\ee
%where $0\leq d_1\leq d_2\leq\cdots\leq d_N$.
Assume that $\gamma^{i}_t$ follows an identical distributed realization
\be
\mathcal{P}(\gamma^{i}_t=0)=p_i,
~\mathcal{P}(\gamma^{i}_t=1)=1-p_i,
~p_i\in[0,1].
\ee
When we define $\omega_t^i=\gamma^{i}_t-(1-p_i)$, WNCS in \eqref{NCS-original} can be viewed as a special case of multi-delay stochastic system \eqref{system01}.
\end{remark}

\begin{remark}
For the multi-delay stochastic systems, \cite{lilin2016,lilin2021} addresses the LQ optimization problem utilizing the stochastic maximum principle. While dealing with the infinite horizon optimization, a necessary and sufficient stabilization condition is derived in terms of CAREs, where $D+1$ positive solutions are to be determined. Generally, it is impossible to verify whether the given result is valid.
Therefore, our purpose is to seeking more concise criteria for stabilizing system (\ref{system01}).
\end{remark}
\section{Main results}

\subsection{Riccati-type Stabilization Criterion}
Motivated by the predictive control technique in \cite{art1982}, we first adopt a feasible reduction strategy to transform system (\ref{system01}) into a delay-free equivalent form.
Based on that, we then develop the Riccati-type stabilization criterion.
To achieve this goal, we need a basic assumption that \emph{$A$ is invertible}.

Define an auxiliary state
\be
\eta_t\hspace{-0.8mm}=\hspace{-0.8mm}x_t
\hspace{-0.8mm}+\hspace{-0.8mm}
\sum_{\tau=1}^D \le[\sum_{j=\tau}^D A^{\tau-j-1} \left(B_j+\omega^j_{t+j-\tau}C_j \right)  \ri]u_{t-\tau}.\label{system02}
\ee
It follows from (\ref{system01}) that
\bee
\eta_{t+1}
%=&Ax_t+\sum_{\tau=0}^D \le(B_\tau +\omega^\tau_t C_\tau \ri)u_{t-\tau} \nonumber\\
%&+\sum_{\tau=1}^D \le[\sum_{j=\tau}^D A^{\tau-j-1} \left(B_j+\omega^j_{t+1+j-\tau}C_j \right)  \ri]u_{t+1-\tau} \no\\
=&Ax_t+\sum_{\tau=1}^D \le[\sum_{j=\tau}^D A^{\tau-j} \left(B_j+\omega^j_{t+j-\tau}C_j \right)  \ri]u_{t-\tau} \no\\
&+\sum_{j=0}^D A^{-j}\left(B_j+\omega_{t+j}^j C_j \right) u_t,
\eee
which leads to the following delay-free expression
\be
\eta_{t+1}=A\eta_t +\sum_{j=0}^D A^{-j}\left(B_j+\omega_{t+j}^j C_j \right) u_t.  \label{system03}
% H u_t+ \sum_{j=0}^D \omega_{t+j}^j A^{-j} C_j  u_t.
\ee
%with
%\be
%H=\sum_{j=0}^D A^{-j}B_j.\label{define-H}
%\ee
For notational convenience, we define
\bse\bee
\bar{\omega}_t\triangleq &\left(\omega_t^0~\omega_{t+1}^1\cdots \omega_{t+D}^D \right)',~t\geq0\\
\bar{\omega}_{-1}\triangleq &\left(\omega_1^0~\omega_{2}^0~\omega_{2}^1\cdots \omega_{D}^{0}\cdots\omega_{D}^{D-1} \right)',
%\mathcal{G}_t=&\sigma\{\bar{\omega}_s,s=-1,0,\cdots,t\}.
\eee\ese
which defines the $\sigma$-algebra as $\mathcal{G}_t=\sigma\{\bar{\omega}_s,s\in\mathbb{N}_{[-1,t]} \}$.
Different from the fact that $x_t\in\mathcal{F}_{t-1}$, the auxiliary state $\eta_t$ is $\mathcal{G}_{t-1}$-measurable, and the pre-defined $\sigma$-algebras satisfy
\be
\mathcal{F}_t \subset \mathcal{G}_t \subset \mathcal{F}_{t+D}.
\ee
Denote $\hat{\eta}_{t|s}={\bf E}\le[\eta_t|\mathcal{F}_{s} \ri]$.
It is of interest to point out that the stabilization of system (\ref{system01}) is equivalent to that of the delay-free system (\ref{system03}) where the control law is designed to be the feedback of $\hat{\eta}_{t|t-1}$ and $u_{t-\tau}$.

\begin{lemma}\label{lemma01}
System (\ref{system01}) is asymptotically mean-square stabilizable if and only if there exists a feedback control law
\be
u_{t}=L_0 \hat{\eta}_{t|t-1} +\sum_{\tau=1}^D L_\tau u_{t-\tau} \in \mathcal {U}_{ad},\label{input01}
\ee
such that the closed-loop auxiliary system in (\ref{system03}) is asymptotically mean-square stable.
\end{lemma}

\textbf{Proof.}
{\em Sufficiency.}
%Suppose there exists a control policy $u_t$ in \eqref{input01} such that system (\ref{system03}) is stable in the mean-square sense, i.e.
%\be
%\lim_{t\rightarrow \infty} {\bf E}\le(u_{t}'u_{t}\ri)=0, \lim_{t\rightarrow \infty} {\bf E}\le(\eta_{t}'\eta_{t}\ri)=0.
%\ee
First, taking the condition expectation on both sides of (\ref{system02}) w.r.t. $\mathcal{F}_{t-1}$ yields
\bee
\hat{\eta}_{t|t-1}=x_t+\sum_{\tau=1}^D \sum_{j=\tau}^D A^{\tau-j-1} B_j u_{t-\tau}.\label{system04}
\eee
The controller in (\ref{input01}) is equivalently expressed by
\beee
u_{t}=L_0 x_t+\sum_{\tau=1}^D \le(\sum_{j=\tau}^D  L_0A^{\tau-j-1} B_j + L_\tau \ri)u_{t-\tau}, \label{input02}
%=&~K_0 x_t+\sum_{\tau=1}^T K_\tau u_{t-\tau}. \no
\eeee
that belongs to the admissible control set $\mathcal {U}_{ad}$ in (\ref{ad-input-set}) with $K_0=L_0$ and $K_\tau=\sum_{j=\tau}^D L_0A^{\tau-j-1} B_j + L_\tau$.

In what follows, let us define
$\tilde{\eta}_{t|t-1}=\eta_t-\hat{\eta}_{t|t-1}.$
In this case, the orthogonality of $\hat{\eta}_{t|t-1}$ and $\tilde{\eta}_{t|t-1}$ can be obtained,
that is
\beee
{\bf E} \Big(\hat{\eta}_{t|t-1}'\tilde{\eta}_{t|t-1} \Big) =& {\bf E} \le[ {\bf E} \le(\hat{\eta}_{t|t-1}' (\eta_t-\hat{\eta}_{t|t-1}) | \mathcal{F}_{t-1} \ri) \ri]=0
%=& {\bf E} \Big(\hat{\eta}_{t|t-1}'\hat{\eta}_{t|t-1}\Big)- {\bf E} \Big(\hat{\eta}_{t|t-1}'\hat{\eta}_{t|t-1} \Big)  =0
%=&0,
\eeee
which leads to
%By utilizing, we obtain
\be
{\bf E}(\eta'_t \eta_t)= {\bf E} (\hat{\eta}_{t|t-1}'\hat{\eta}_{t|t-1})+{\bf E}(\tilde{\eta}_{t|t-1}' \tilde{\eta}_{t|t-1}).
\ee
Hence,
%Define $\bar{M}_{\tau}=\sum_{j=\tau}^T A^{\tau-j-1}B_j$ and $\bar{N}_\tau=\sum_{j=\tau}^T \omega^j_{t+j-\tau}\\A^{\tau-j-1}C_j$.
by (\ref{system04}), one obtains that
%\be
%x_t=\eta_t-\sum_{\tau=1}^D \sum_{j=\tau}^D A^{\tau-j-1} B_j u_{t-\tau}
%%\sum_{\tau=1}^T \le[ \sum_{j=\tau}^T A^{\tau-j-1} \left(B_j+\omega^j_{t+j-\tau}C_j \right) \ri]u_{t-\tau},
%%\sum_{\tau=1}^T \le(\bar{M}_{\tau}+ \bar{N}_\tau \ri) u_{t-\tau},
%\ee
%and
\bee\label{inequ}
&{\bf E}\le(x_{t}'x_{t}\ri) \leq
{\bf E}
\Bigg\{
\le(\eta_{t}'\eta_{t}\ri)^{\frac{1}{2}}
+\sum_{\tau=1}^D \Big[
u_{t-\tau}' \le( \sum_{j=\tau}^D A^{\tau-j-1} B_j \ri)' \no\\
&\times
\le( \sum_{j=\tau}^D A^{\tau-j-1} B_j\ri) u_{t-\tau} \Big]^{\frac{1}{2}} \Bigg\} ^2 \no\\
&\leq (D+1)\Bigg\{ {\bf E}\le(\eta_{t}'\eta_{t}\ri)
+\sum_{\tau=1}^D {\bf E} \Big[u_{t-\tau}' \le( \sum_{j=\tau}^D A^{\tau-j-1}B_j \ri)'
\no\\
&\times
\le( \sum_{j=\tau}^D A^{\tau-j-1}B_j \ri) u_{t-\tau}\Big] \Bigg\}.
%&+\sum_{\tau=1}^T {\bf E} \Big[u_{t-\tau}' \left(\sum_{j=\tau}^T \sigma_j^2 C_j'(A')^{\tau-j-1} A^{\tau-j-1}C_j \right) u_{t-\tau}  \Big]
\eee
%\bee
%&{\bf E}\le(x_{t}'x_{t}\ri) \no\\
%\leq &
%{\bf E}
%\Bigg\{
%\le(\eta_{t}'\eta_{t}\ri)^{\frac{1}{2}}
%+\sum_{\tau=1}^T \Big[
%u_{t-\tau}' \le( \sum_{j=\tau}^T A^{\tau-j-1} \left(B_j+\omega^j_{t+j-\tau}C_j \right) \ri)'  \no\\
%&\times  \le( \sum_{j=\tau}^T A^{\tau-j-1} \left(B_j+\omega^j_{t+j-\tau}C_j \right) \ri) u_{t-\tau} \Big]^{\frac{1}{2}} \Bigg\} ^2 \no\\
%\leq & (T+1)\Bigg\{ {\bf E}\le(\eta_{t}'\eta_{t}\ri)
%+\sum_{\tau=1}^T {\bf E} \Big[u_{t-\tau}' \le( \sum_{j=\tau}^T A^{\tau-j-1}B_j \ri)' \no\\
%&\times \le( \sum_{j=\tau}^T A^{\tau-j-1}B_j \ri) u_{t-\tau}\Big]  \no\\
%&+\sum_{\tau=1}^T {\bf E} \Big[u_{t-\tau}' \left(\sum_{j=\tau}^T \sigma_j^2 C_j'(A')^{\tau-j-1} A^{\tau-j-1}C_j \right) u_{t-\tau}  \Big] \Bigg\}.
%\eee
Since $\le( \sum_{j=\tau}^D A^{\tau-j-1}B_j \ri)'\le( \sum_{j=\tau}^D A^{\tau-j-1}B_j \ri)\geq 0$,
%and $C_j'(A')^{\tau-j-1} A^{\tau-j-1}C_j\geq 0$,
there exists a positive scalar $\gamma>0$ such that
\bee
{\bf E}\le(x_{t}'x_{t}\ri) \leq (D+1)\Bigg\{ {\bf E}\le(\eta_{t}'\eta_{t}\ri)
+\gamma\sum_{\tau=1}^D {\bf E} \left(u_{t-\tau}'u_{t-\tau}\right) \Bigg\}, \no
\eee
which implies that $\lim_{t\rightarrow \infty} {\bf E}\le(x_{t}'x_{t}\ri)=0$.

{\em Necessity.} This part can be derived from \eqref{inequ}.
\ok

Based on Lemma \ref{lemma01}, our objective in Problem 1 can be reformulated as seeking an admissible control law to stabilize auxiliary system \eqref{system03}.
Below, we propose the Riccati-type stabilization criterion.

\begin{theorem}\label{riccati-theorem}
System (\ref{system01}) is asymptotically mean-square stabilizable if and only if for any $Q>0$ and $R>0$, there exists a unique positive definite solution $Z>0$ satisfying the following nonlinear DDARE
\be
-Z+A'{Z}A+{Q}-A'ZL\Psi^{-1}L'{Z}A=0,\label{DARE-trans}
\ee
where
\bse
\bee
\Psi=& L'{Z}L+\sum_{\tau=0}^D \sigma_\tau^2 C_\tau' (A')^D {Z} A^D C_\tau+U_{{R},{Q}},  \label{DARE-trans-parameter-P} \\
L=&\sum_{j=0}^D A^{D-j}B_j, \label{DARE-trans-parameter-L}\\
U_{{R},{Q}}=&{R}+\sum_{\tau=1}^D \sum_{h=1}^\tau\sigma_\tau^2 C'_\tau (A')^{D-h} {Q} A^{D-h} C_\tau. \label{DARE-trans-parameter-V}
\eee
\ese
Moreover, the stabilizing control policy is given as
\bee
\hspace{-2.5mm}
u^*_t
%=&-\Psi^{-1}L'ZA^{D+1} \hat{\eta}_{t|t-1} \no\\
=-\Psi^{-1}L'ZA^{D} \left( \hspace{-1mm} Ax_t+\sum_{\tau=1}^D \sum_{j=\tau}^D A^{\tau-j} B_j u^*_{t-\tau}\right).
\eee
\end{theorem}

\textbf{Proof.}
See Appendix  \ref{app02}.
\ok

\begin{remark}
In Theorem 1, we have proposed the DDARE-based stabilization condition, which is necessary and sufficient.
However, since the developed Riccati equation is nonlinear, how to calculate the unique positive solution is challenging and to be solved.
\end{remark}

\subsection{Lyapunov-type Stabilization condition}
In this subsection, we propose a Lyapunov-type necessary and sufficient stabilization condition, which can be verified availably by the LMI feasibility test.
We reveal an interesting fact that the stabilization of system \eqref{system01} is equivalent to stabilizing another delay-free model of same dimensions
\be
\beta_{t+1}=A\beta_t+Lv_t+
\sum_{\tau=0}^D \omega^\tau_t A^DC_\tau v_{t},\label{system05}
\ee
where $L$ is defined in \eqref{DARE-trans-parameter-L}.
%the input control $v_t$ is designed to be the feedback of the state information, i.e..
We remark that system \eqref{system05} is said to be asymptotically mean-square stabilizable, if there exists a state feedback law $v_t=K\beta_t$ such that the following closed-loop system is stable \cite{zhouxy2000},
\be
\beta_{t+1}=(A+LK)\beta_t+
\sum_{\tau=0}^D \omega^\tau_t A^DC_\tau K\beta_t.\label{system05-close}
\ee

\begin{theorem}\label{Lyapunov-theorem}
The following statements are equivalent.\\
$a)$ System \eqref{system01} is asymptotically mean-square stabilizable. \\
$b)$ System \eqref{system05} is asymptotically mean-square stabilizable with $v_t=K\beta_t$. \\
$c)$ For any $Q>0$, there exist matrices $K$ and $P>0$ satisfying the following delay-dependent Lyapunov equation
\beee
P=Q+(A+LK)'P(A+LK)
+\sum_{\tau=0}^D \sigma^2_\tau K'C_\tau' (A')^D PA^DC_\tau K.
\eeee
$d)$ There exist matrices $K$ and $P>0$ satisfying the following delay-dependent Lyapunov inequality
\beee
P>& (A+LK)'P(A+LK)
+\sum_{\tau=0}^D \sigma^2_\tau K'C_\tau' (A')^D PA^DC_\tau K.\label{inequality}
\eeee
\end{theorem}

\textbf{Proof.}
Similar to Theorem 1 in \cite{zhouxy2000}, we have that $b) \Leftrightarrow c)$ and $c) \Leftrightarrow d)$.
Based on Theorem \ref{riccati-theorem}, we only need to prove that the stabilization of system \eqref{system05} is equivalent to DDARE \eqref{DARE-trans} has a unique positive solution.

$a)\Rightarrow b)$
Suppose system \eqref{system01} is stabilizable.
Let us define
\bee
V_t(\beta)= {\bf E} \le(\beta_t'Z\beta_t \ri),
\eee
where $Z>0$ is the unique positive solution of DDARE \eqref{DARE-trans}.
It follows that
\bee
&V_t(\beta)-V_{t+1}(\beta) ={\bf E} \Big\{\beta_t' Q \beta_t+v_t' U_{R,Q} v_t \no\\
& -\big( v_t+\Psi^{-1}L'PA {\beta}_{t} \big)'\Lambda\big( v_t+\Psi^{-1}L'PA {\beta}_{t}  \big) \Big\}.
\eee
%where
%\bse
%\bee
%Q_1=&\hat{Q}\succ0, \\
%R_1=&\sum_{\tau=1}^D \sum_{h=1}^\tau\sigma_\tau^2 C_\tau (A')^{D-h} \hat{Q} A^{D-h} C_\tau+\hat{R}\succ0.
%\eee
%\ese
Utilizing the control law $v_t=-\Psi^{-1}L'PA {\beta}_{t}$ yields that $V_t(\beta)-V_{t+1}(\beta)>0$, which indicates that system \eqref{system05} is stabilizable.

$b)\Rightarrow a)$
%To begin with,
Define the following finite horizon index function
\bee
J_T(\beta_0,v_t)
=\sum_{t=0}^T {\bf E}\le[\beta_t'Q\beta_t+v_t'U_{R,Q}v_t\ri].
\eee
By utilizing the matrix version of the maximum principle in \cite{zhouxy2002}, one obtains that the finite horizon LQ problem admits a unique optimal control
\be
v_t=-\Psi_{t+1}^{-1}L'P^T_{t+1}A\beta_t,
\ee
where the unique solution $P^T_t$ solves the following generalized difference Riccati equation (GDRE)
\bse
\bee
P^T_t=&A'P^T_{t+1} A+Q-A'P^T_{t+1} L\Psi_{t+1}^{-1}L'P^T_{t+1}A,\label{GDRE} \no\\
\Psi_{t+1}=& L'P^T_{t+1}L+\sum_{\tau=0}^D \sigma_\tau^2 C_\tau' (A')^DP^T_{t+1} A^D C_\tau+U_{R,Q}, \no
\eee
\ese
with the terminal condition $P^T_{T+1}=0$.
Similar to Theorem 1 in \cite{zhangwh2008}, one obtains that the stabilization of system \eqref{system05} guarantee the existence of the limit $P=\lim_{T\rightarrow\infty}P^T_t$.
Moreover, by taking the limit on both sides of \eqref{GDRE}, one obtains that $P>0$ is the unique positive solution satisfying the following GARE
\bee
P=A'P A+Q-A'P L\Psi^{-1}L'PA,
\eee
which is equivalent to DDARE \eqref{DARE-trans} with $Z=P$.
\ok

Utilizing Schur complement decomposition and matrix transformation technique, we give the following LMI-based criterion directly, which can be verified by some LMI solvers.

\begin{corollary}
System \eqref{system01} is stabilizable if and only if there exist matrices $Y$ and $S<0$ satisfying
\bee\label{lmi}
\left[
  \begin{array}{ccccc}
    S & * & * & \cdots & * \\
    AS+LY & S & * & \cdots & * \\
    \sigma_0 A^DC_0 Y & 0 & S & \cdots & * \\
    \vdots & \vdots & \vdots & \ddots & \vdots \\
    \sigma_D A^DC_D Y & 0 & 0 & \cdots & S \\
  \end{array}
\right]< 0,
\eee
where $Y=KP^{-1}$, $S=-P^{-1}<0$ and
$*$ represents the corresponding transpose part.
\end{corollary}

\begin{remark}
Specifically, let us set $B_i=C_i=0$, $i\in\mathbb{N}_{[0,D-1]}$. System \eqref{system01} is reduced to a single-delay stochastic system with control-dependent noise
\be
x_{t+1}=Ax_t+\le(B_D+\omega_t C_D \ri)u_{t-D}.\label{system01re}
\ee
Based on Theorem \ref{Lyapunov-theorem}, system \eqref{system01re} is stabilizable if and only if there exists a unique ${P}>0$ satisfying
\bee
P=Q+(A+LK)'P(A+LK)
+\sigma^2_D K'C_D' (A')^D PA^DC_D K,
\eee
which corroborates with the developed Lyapunov-type criterion in Theorem 2 in \cite{tan2019}.
\end{remark}

\section{Application}

This section focuses on the delay margin in the following restricted structure
\be
x_{t+1}=Ax_t+\omega_t^0 C_0 u_{t}+ \le(\bar{B} +\omega_t^D\bar{C} \ri)u_{t-D}, \label{system01single}
\ee
which will be described as $[A,C_0;\bar{B},\bar{C}]_D$ hereinafter.
The problem with our concern is to determine the largest delay range within which $[A,C_0;\bar{B},\bar{C}]_D$ is stabilizable.

Since for a stable system, an arbitrarily large delay margin can be achieved by zero input.
To proceed we state the following assumptions.

\noindent{\bf Assumption 1} $A$ is invertible and unstable. %and $\bar{B}$ has full-column rank.

\noindent{\bf Assumption 2} System $[A,C_0;\bar{B},\bar{C}]_0$ is stabilizable and the following delay-free system $(A,C_0)$ is destabilizable \cite{tan2022},
\be
x_{t+1}=Ax_t+\omega_t^0 C_0 u_{t}.
\ee
It is remarkable that Assumption 2 indicates that system $[A,C_0;\bar{B},\bar{C}]_D$ is stabilizable in delay-free case and it can not be stabilized when $D$ tends to infinity; See Fig. 2.

\begin{figure}[thpb]
      \centering
      \includegraphics[scale=0.38]{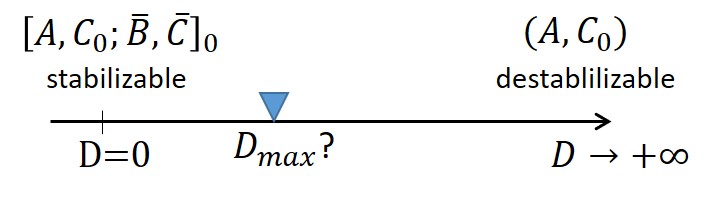}
      \centering\caption{\scriptsize{Stabilization classification for system $[A,C_0;\bar{B},\bar{C}]_D$}}
\end{figure}

For the sake of discussion, we introduce the following nonlinear Riccati operator $\mathcal {R}^R_{D}(\cdot)$ as
\be
\mathcal {R}^R_{D}(Z)\triangleq  A'{Z}A+{Q}-A'Z\bar{B}\Psi^R_{D}(Z)^{-1}\bar{B}'{Z}A,\label{DARE-trans-operator}
\ee
where
\beee
\Psi^R_{D}(Z) =&\bar{B}'{Z}\bar{B}+\sigma_0^2 C_0' (A')^D {Z} A^D C_0+\sigma_D^2 \bar{C}' (A')^D {Z} A^D \bar{C}\no\\
&+\sum_{h=1}^D\sigma_D^2 \bar{C}' (A')^{D-h} {Q} A^{D-h} \bar{C}+R. \label{DARE-trans-parameter-P-up}
\eeee
Moreover, define an auxiliary Lyapunov operator $\mathcal {J}_K(\cdot)$
\beee
&\mathcal {J}_K(Z) \triangleq
(A+LK)'Z(A+LK) +\sum_{\tau=0}^D \sigma^2_\tau  K'C_\tau'(A')^D Z\no\\
&\times A^DC_\tau K+Q+K'(R+\sum_{h=1}^D\sigma_D^2 \bar{C}' (A')^{D-h} {Q} A^{D-h})K . %\label{lyapunov-J}
\eeee
Specifically, denote $K_Z=-\Psi^R_{D}(Z)^{-1}\bar{B}'{Z}A$, and the corresponding DDARE for $[A,C;\bar{B},\bar{C}]_D$ can be expressed as
\be
Z=\mathcal {R}^R_{D}(Z)=  \mathcal {J}_{K_Z}(Z). \label{lyapunov-J-2}
\ee
%and the stabilizing control policy is given as
%\bee\label{stable-input-02}
%u^*_t
%=-\Psi^{-1}\bar{B}'Z\left( A^{D+1}x_t+\sum_{\tau=1}^D A^{\tau} \bar{B} u^*_{t-\tau}\right).
%\eee
In addition, for any given $K$, using the completing square approach yields that
\bee
\mathcal {J}_{K}(Z)=&A'{Z}A+{Q}-A'ZL\Psi^{-1}L'{Z}A \no\\
&+(K-K_Z)'\Psi (K-K_Z)   \no\\
&\geq \mathcal {J}_{K_Z}(Z)=\mathcal {R}^R_{D}(Z).\label{lyapunov-J-3}
\eee

\begin{theorem}\label{existence-theorem}
Under Assumptions 1-2, there exists a unique delay margin $D_{\max}\geq 0$ such that $[A,C_0;\bar{B},\bar{C}]_D$ is asymptotically mean-square stabilizable for $D\in \mathbb{N}_{[0,D_{\max}]}$, and it can not be stabilized for each $D>D_{\max}$.
\end{theorem}

\textbf{Proof.}
To proceed we first show that the stabilization of system $[A,C_0;\bar{B},\bar{C}]_{D+1}$ can guarantee that of $[A,C_0;\bar{B},\bar{C}]_{D}$.
Select some positive definite matrices $R$ and $Q$ such that
\be
\hat{R}\triangleq R-\sigma_0^2 C_0' (A')^D Q A^D C_0>0.
\ee
For such $R>0$ and $Q>0$, there exists a unique positive solution $Z$ satisfying $Z=\mathcal {R}^R_{D+1}(Z)$.
It follows that
$Z\leq A'ZA+{Q}$,
which implies that
\bee
\sigma_0^2 C_0' (A')^D Z\ A^D C_0 -\sigma_0^2 C_0' (A')^D Q A^D C_0 \no\\
\leq
\sigma_0^2 C_0' (A')^{D+1} Z\ A^{D+1} C_0, \\
\sigma_D^2 \bar{C}' (A')^D Z\ A^D \bar{C} -\sigma_D^2 \bar{C}' (A')^D Q A^D \bar{C} \no\\
\leq
\sigma_D^2 \bar{C}' (A')^{D+1} Z\ A^{D+1}\bar{C}.
\eee
We further derive that $\Psi^R_{D+1}(Z)\geq \Psi^{\hat{R}}_{D}(Z).$
%\bee
%&\Psi^R_{D+1}(Z)  \no\\
%=& \bar{B}'Z\bar{B}+ \sigma_0^2 C_0' (A')^{D+1} Z A^{D+1} C_0+\sigma_D^2 \bar{C}' (A')^{D+1} Z \no\\
%&\times A^{D+1} \bar{C}+R+ \sum_{h=1}^{D+1}\sigma_D^2 \bar{C}' (A')^{D+1-h} {Q} A^{D+1-h} \bar{C} \no\\
%\geq &
%\bar{B}'Z\bar{B}+ \sigma_0^2 C_0' (A')^D Z A^D C_0 +\sigma_D^2 \bar{C}' (A')^D Z A^D \bar{C}+\hat{R}\no\\
%& +\sum_{h=1}^D\sigma_D^2 \bar{C}' (A')^{D-h} {Q} A^{D-h} \bar{C}
%=\Psi^{\hat{R}}_{D}(Z).
%\eee
Then, we can get that
$Z=\mathcal {R}^R_{D+1}(Z)\geq \mathcal {R}^{\hat{R}}_{D}(Z)$.
If we define $Z_{k+1}=\mathcal {R}^{\hat{R}}_D(Z_k)$ with $Z_0=Z$, it follows that
\be
Z_0\geq\mathcal {R}^{\hat{R}}_D(Z_0)=Z_1.
\ee
With the help of the pre-defined operator $\mathcal {J}_K(\cdot)$ and $K_{Z_k}=-\Psi^{\hat{R}}_{D}(Z_k)^{-1}\bar{B}'{Z_k}A$, we obtain that
\be
Z_{k+1}=\mathcal {R}^{\hat{R}}_{D}(Z_k)=\mathcal {J}_{K_{Z_k}}(Z_k)\geq Q>0.\label{iterative}
\ee
Then, it follows from \eqref{lyapunov-J-3} that
\bee
Z_1=\mathcal {J}_{K_{Z_0}}(Z_0)\geq\mathcal {J}_{K_{Z_0}}(Z_1)\geq\mathcal {J}_{K_{Z_1}}(Z_1)=Z_2,
\eee
%\bee
%Z_1=\mathcal {R}^{\hat{R}}_D(Z_0)=\mathcal {J}_{K_{Z_0}}(Z_0)
%\geq\mathcal {J}_{K_{Z_0}}(Z_1).
%\eee
%\bee
%Z_1\geq\mathcal {J}_{K_{Z_0}}(Z_1)\geq\mathcal {J}_{K_{Z_1}}(Z_1)=\mathcal {R}^{\hat{R}}_D(Z_1)=Z_2,
%\eee
which results in $Z_k\geq Z_{k+1}$ by induction.
i.e., $\{Z_k\}$ is a monotone decreasing and bounded sequence.
Therefore, there exists a $\hat{Z}=\lim_{k\rightarrow\infty}Z_k$ satisfying
\be
\hat{Z}=\lim_{k\rightarrow\infty}\mathcal {R}^R_{D}(Z_k)=\mathcal {R}^R_{D}(\hat{Z})
=\mathcal {J}_{K_{\hat{Z}}}(\hat{Z})>0,
\ee
which means that system $[A,C_0;\bar{B},\bar{C}]_{D}$ is stabilizable.
%based on Theorem \ref{riccati-theorem}.
% Corollary \ref{riccati-corollary}.

When $D=0$, system $[A,C_0;\bar{B},\bar{C}]_{0}$ is stabilizable on account of Assumption 2.
And when $D$ tends to $+\infty$, system $[A,C_0;\bar{B},\bar{C}]_{D}$ is equivalent to delay-free system $(A,C_0)$, that cannot be stabilized.
To sum up, the existence and uniqueness of delay margin $D_{\max}\geq 0$ is proposed.
\ok

%\begin{proof}
%See Appendix \ref{app03}.
%\end{proof}

\begin{remark}\label{remark-c}
The proof in Theorem \ref{existence-theorem} produces a numerically iterative algorithm for computing the definite positive solution to DDARE \eqref{lyapunov-J-2}.
Suppose system $[A,C_0;\bar{B},\bar{C}]_D$ is stabilizable.
For any initial value $Z_0\geq0$ satisfying $Z_0\geq\mathcal {R}^{{R}}_D(Z_0)$ with $Q>0$ and $R>0$, the solution satisfying the recurrence formula $Z_{k+1}=\mathcal {R}^{{R}}_D(Z_k)$ converges to the unique positive solution in \eqref{lyapunov-J-2}.
However, for general nonlinear DDARE \eqref{DARE-trans}, how to compute the value of the unique positive solution is
still unknown, which defines a challenging work direction.
\end{remark}

Below we study an uncoupled system, where
\bee
A=diag\{a_1,\cdots,a_n\},~C_0=diag\{c_1,\cdots,c_n\}, \no\\
\bar{B}=diag\{\bar{b}_1,\cdots,\bar{b}_n\},~\bar{C}=diag\{\bar{c}_1,\cdots,\bar{c}_n\}. \no
\eee
Hence, system $[A,C_0;\bar{B},\bar{C}]_{D}$ is stabilizable if and only if each scalar sub-system $[a_i,c_i;\bar{b}_i,\bar{c}_i]_{D}$ is stabilizable, and $D_{\max}$ satisfies
\be
D_{\max}=\min_{1\leq i\leq n}D^i_{\max},
\ee
where $D^i_{\max}$ is the delay margin for $[a_i,c_i;\bar{b}_i,\bar{c}_i]_{D}$.

Note that when $a_i^2<1$, the scalar sub-system $[a_i,c_i;\bar{b}_i,\bar{c}_i]_{D}$ is stabilizable for any bounded $D>0$.
To obtain a bounded delay margin, we assume that $A$ has at least one unstable eigenvalue $a^2_i\geq 1$ and $\bar{b}_j\neq 0$, $j\in\mathbb{N}_{[1,n]}$.
Denote $U_A$ to be the unstable eigenvalue set of the diagonal matrix $A$, i.e.,
\be
U_A\triangleq \{a_i:a^2_i\geq 1,i\in\mathbb{N}_{[1,n]}\}.
\ee
We are in a position to derive the exact delay margin.
%in terms of the system matrices and the probability statistics of multiplicative noises.

\begin{theorem}\label{delay-margin-theorem}
System $[A,C_0;\bar{B},\bar{C}]_{D}$ is stabilizable if and only if $0\leq D<D_{\max}$, where $D_{\max}$ is the  delay margin and satisfies the following conditions:
\\
a) If for each $a_i\in U_A$, $a_i^2=1$, then $D_{\max}=+\infty$.
\\
b) If for each $a_i\in U_A$, $a_i^2<\frac{\bar{b}^2_i}{\sigma_0^2 c^2_i+\sigma_D^2 \bar{c}^2_i}+1$,
and there exists a $a_j\in U_A$ such that $a_j^2>1$,
%there exists a scalar $a_j\in U(A)$ such that $1<a^2_j<{\rm SNR}_j+1$,
then
\beee
D_{\max}
=\min_{a_i\in U_A,a_i^2>1}\frac{\ln(\bar{b}^2_i)-\ln(\sigma_0^2 c^2_i+\sigma_D^2 \bar{c}^2_i)-\ln(a_i^2-1)}{\ln(a_i^2)}.
\eeee
c) If there exists a $a_i\in U_A$ satisfying $a_i^2\geq \frac{\bar{b}^2_i}{\sigma_0^2 c^2_i+\sigma_D^2 \bar{c}^2_i}+1$, then
$D_{\max}=0$.
\end{theorem}

\textbf{Proof.}
See Appendix  \ref{app03}.
\ok

\section{Simulation}

In this section, we give two examples to demonstrate our theoretical analysis.

\begin{example}
For simplicity, we consider the following system \eqref{system01} with $D=2$ and
\beee
&A=\left[ {\begin{array}{*{20}{c}}
1& 2&3\\
0&2 &2 \\
0&0&1
\end{array}} \right],
B_0=\left[ {\begin{array}{*{20}{c}}
2 &3  \\
1 &1  \\
1 &4
\end{array}} \right], \\
&B_1=\left[ {\begin{array}{*{20}{c}}
2 &4 \\
 2 &5 \\
 2 &2
\end{array}} \right],
B_2=\left[ {\begin{array}{*{20}{c}}
3 &4 \\
1 &3 \\
3 &5
\end{array}} \right], \\
&C_0=\left[ {\begin{array}{*{20}{c}}
 5     &5      \\
 40    &-3    \\
  3     &2
\end{array}} \right],
C_1=\left[ {\begin{array}{*{20}{c}}
2     &5      \\
2     &0      \\
4     &4
\end{array}} \right],
C_2=\left[ {\begin{array}{*{20}{c}}
0     &3 \\
2     &0 \\
1     &2 \\
\end{array}} \right],
\eeee
Assume $\sigma_\tau^2={\bf E}\le(\omega_t^\tau \omega_t^\tau\ri)=1,~\tau=0,1,2$,
and the initial condition is $x_0=[5~-3~10]'$ and $u_{-i}=[0~0]',~i=1,2$.
Utilizing the developed LMI in Corollary 1 yields the stabilizing feedback gain is
\beee
K=\left[ {\begin{array}{*{20}{c}}
   -0.0001 &   -0.0028  & -0.0105 \\
    0.0004  &  0.0132  &  0.0571
\end{array}} \right].
\eeee
In this case, the control policy $u_t=K\hat{\eta}_{t|t-1}$ can stabilize system \eqref{system01} as shown in Fig. 3.\ok

\begin{figure}[thpb]
      \centering
      \includegraphics[scale=0.3]{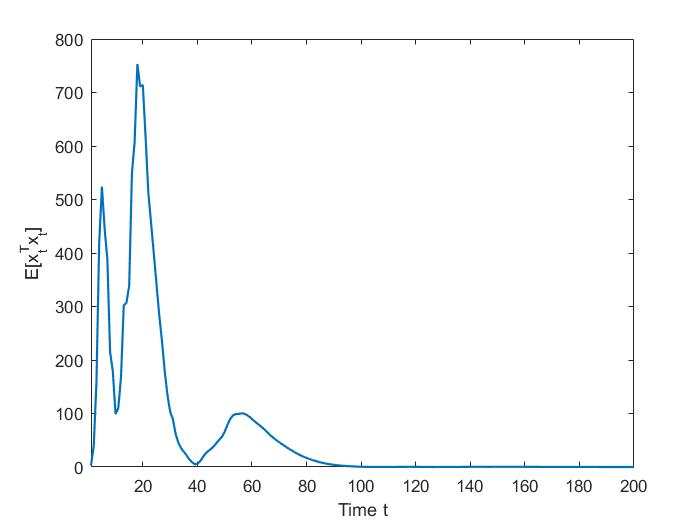}
      \centering\caption{\scriptsize{State response of system (1) in Example 1}}
\end{figure}

\end{example}

\begin{example}
Consider system $[A,C_0;\bar{B},\bar{C}]_2$ with the following parameters
\beee
A=\left[ {\begin{array}{*{20}{c}}
1.1& 0\\
0&1.2
\end{array}} \right],
C_0=\left[ {\begin{array}{*{20}{c}}
1& -0.5\\
0& 1.25
\end{array}} \right], \\
\bar{B}=\left[ {\begin{array}{*{20}{c}}
2& 1\\
1&-1
\end{array}} \right],
\bar{C}=\left[ {\begin{array}{*{20}{c}}
2& 0\\
0&3
\end{array}} \right].
\eeee
Assume $\sigma_\tau^2={\bf E}\le(\omega_t^\tau \omega_t^\tau\ri)=1,~\tau=0,2$,
and the initial condition is $x_0=[1~1]'$ and $u_{-i}=[0~0]',~i=1,2$.
Utilizing the iterative algorithm proposed in Remark \ref{remark-c} yields that the unique positive solution to DDARE \eqref{lyapunov-J-2} is
\beee
Z=
\left[ {\begin{array}{*{20}{c}}
82.7362& -257.7524\\
 -257.7524 & 859.1263
\end{array}} \right]>0.
\eeee
In this case, the stabilizing control policy in Theorem \ref{riccati-theorem} is
%can be reduced to
\beee
u^*_t
=-\Psi^{-1}\bar{B}'Z\left( A^{D+1}x_t+\sum_{\tau=1}^D A^{\tau} \bar{B} u^*_{t-\tau}\right),
\eeee
which stabilizes the considered system as shown in Fig. 4.\ok
\begin{figure}[thpb]
      \centering
      \includegraphics[scale=0.3]{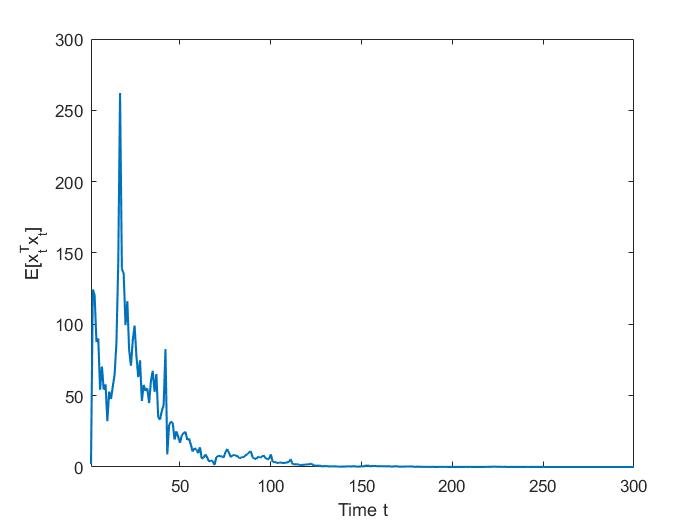}
      \centering\caption{\scriptsize{State response of system $[A,C_0;\bar{B},\bar{C}]_2$ in Example 2}}
\end{figure}

\end{example}

\section{Conclusion}
In this paper, we are concerned with the stabilization problem of multi-delay stochastic system with multiplicative noises in control variables.
By utilizing a novel reduction method, we derive a set of necessary and sufficient conditions for stabilizing such systems in terms of DDARE and DDLE.
Under some basic assumptions, we prove that the unique delay margin exists, and further propose its explicit computation formula for uncoupled system.
%A promising future direction is applying the developed method to deal with continuous-time multi-delay dynamics.

%%%%%%%%%%%%%%%%%%%%%%%%%%%%%%%%%%%%%%%%%%%%%%%%%%%%%%%%%%%%%%%%%%%%%%%%%%%%%%%%%%%%%%%%%%%%%%%%%%%%%%%%%%%%%%%%%%%%%%%%%%%%%%%%%%%%%%%%%%%

%\section*{Appendix}
\appendices
%\section{Proof of Lemma \ref{lemma01}}\label{app01}

\section{Proof of Theorem \ref{riccati-theorem}}\label{app02}
\textbf{Proof.}
First, we give an equivalent form of DDARE \eqref{DARE-trans}.
Define $P=(A')^{D}ZA^{D}$ and $\hat{Q}=(A')^D{Q}A^D$.
Pre-multiplying $(A')^{D}$ and post-multiplying $A^{D}$ on both sides of DDARE \eqref{DARE-trans}, we obtain that
\bee
0=&(A')^{D}(-Z+A'{Z}A+Q)A^{D}-A'(A')^{D}ZA^{D} \no\\
&\times\left(\sum_{j=0}^D A^{-j}B_j\right)\Psi^{-1} \left(\sum_{j=0}^D A^{-j}B_j\right)(A')^{D}ZA^{D}A \no\\
=&-P+A'PA+\hat{Q}-A'PH\Lambda^{-1}H'PA, \label{DARE}
\eee
where
\bse
\bee
\Lambda=& H'PH+\sum_{\tau=0}^D \sigma_\tau^2 C_\tau'PC_\tau+W_{\hat{R},\hat{Q}}, \label{dare-parameter-L} \\
H=&\sum_{j=0}^D A^{-j}B_j, \label{dare-parameter-H} \\
W_{R,\hat{Q}}=&R+\sum_{\tau=1}^D \sum_{h=1}^\tau\sigma_\tau^2 C_\tau (A')^{-h}\hat{Q}A^{-h} C_\tau.  \label{dare-parameter-W}
\eee
\ese
Below, we prove that system \eqref{system03} is stabilizable if and only if DDARE \eqref{DARE} has a unique positive definite solution.

{\em Sufficiency.}
Due to the equivalent proposition between DDAREs \eqref{DARE-trans} and \eqref{DARE}, we are in a position to construct a Lyapunov function to guarantee $\lim_{t\rightarrow \infty} {\bf E}\le(\eta_{t}'\eta_{t}\ri)=0$
with $u_t=-\Lambda^{-1}H'PA \hat{\eta}_{t|t-1}$.

Step 1:
Based on DDARE \eqref{DARE}, define the following delay-dependent Lyapunov function
\bee
&V_t(\eta)={\bf E} \bigg\{\eta_t' \Big[ (A')^D PA^D+\sum_{\tau=0}^{D-1} (A')^\tau \hat{Q}A^\tau \Big] \eta_t \no\\
&-\sum_{\tau=0}^{D-1} \hat{\eta}_{t|t+\tau-1}' (A')^{\tau+1}PH \Lambda^{-1} H'PA^{\tau+1} \hat{\eta}_{t|t+\tau-1}
\bigg\}.
\eee
The orthogonality of $\hat{\eta}_{t|t+\tau-1}$ and $\tilde{\eta}_{t|t+\tau-1}$ yields that
\beee
&V_t(\eta)
%=& {\bf E} \bigg\{\eta_t' \Big[ (A')^D PA^D+\sum_{\tau=0}^{D-1} ( (A')^\tau \hat{Q}A^\tau- (A')^{\tau+1} \no\\
%&\times PH \Lambda^{-1} H'PA^{\tau+1} )\Big] \eta_t+\sum_{\tau=0}^{D-1} \tilde{\eta}_{t|t+\tau-1}' \no\\
%& \times(A')^{\tau+1}PH \Lambda^{-1} H'PA^{\tau+1} \tilde{\eta}_{t|t+\tau-1}
%\bigg\} \no\\
= {\bf E} \bigg\{\eta_t' \Big[ (A')^D PA^D+\sum_{\tau=0}^{D-1} ( (A')^\tau \hat{Q}A^\tau + (A')^{\tau}
\no\\
&\times
\le(P-A'PA-\hat{Q} \ri)  A^{\tau} )\Big] \eta_t \no \\
&+\sum_{\tau=0}^{D-1} \tilde{\eta}_{t|t+\tau-1}'  (A')^{\tau+1}PH \Lambda^{-1} H'PA^{\tau+1} \tilde{\eta}_{t|t+\tau-1}
\bigg\} \no\\
=& {\bf E} \bigg\{\eta_t'P\eta_t+\sum_{\tau=0}^{D-1} \tilde{\eta}_{t|t+\tau-1}' (A')^{\tau+1}PH \Lambda^{-1} H'PA^{\tau+1} \tilde{\eta}_{t|t+\tau-1}
\bigg\}.
\eeee
Because of $P>0$ and $(A')^{\tau+1}PH \Lambda^{-1} H'PA^{\tau+1}\geq 0$, it follows that $V_t(\eta)>0$ for any $\eta_t\neq 0$.

Step 2: Now move to calculate the difference of $V_{t}(\eta)$ and $V_{t+1}(\eta)$ as follows
\beee\label{lyapunov02}
&V_{t}(\eta)-V_{t+1}(\eta) \no\\
=& {\bf E} \Bigg\{ \eta_t' \bigg[ (A')^D PA^D+\sum_{\tau=0}^{D-1} (A')^\tau \hat{Q}A^\tau \bigg] \eta_t
\no\\
&
-\sum_{\tau=0}^{D-1} \hat{\eta}_{t|t+\tau-1}' (A')^{\tau+1}PH \Lambda^{-1} H'PA^{\tau+1} \hat{\eta}_{t|t+\tau-1} \no\\
&+\sum_{\tau=0}^{D-1} \hat{\eta}_{t+1|t+\tau}' (A')^{\tau+1}PH \Lambda^{-1} H'PA^{\tau+1} \hat{\eta}_{t+1|t+\tau}
\no\\
&-\eta_{t+1}' \bigg[ (A')^DPA^D+\sum_{\tau=0}^{D-1} (A')^\tau \hat{Q}A^\tau \bigg] \eta_{t+1} \Bigg\}.
\eeee
Since $\eta_t$ is $\mathcal{G}_{t-1}$-measurable and $\mathcal{G}_{t-1} \subset \mathcal{F}_{t+D-1}$, we have
$\hat{\eta}_{t|t+D-1}={\bf E} [\eta_t|\mathcal{F}_{t+D-1}]=\eta_t.$
%Substituting (\ref{system02}) into (\ref{lyapunov02}) yields that
%\beee
%&V_{t}(\eta)-V_{t+1}(\eta) \no\\
%=& {\bf E} \Bigg\{  \eta_t' \bigg[ (A')^DPA^D+\sum_{\tau=0}^{D-1} (A')^\tau \hat{Q} A^\tau-(A')^{D+1}PA^{D+1} \no\\
%&- \sum_{\tau=0}^{D-1} (A')^{\tau+1} \hat{Q}A^{\tau+1} \bigg]   \eta_t
%-u_t'  \bigg[ H'(A')^D P A^DH   \no\\
%&+\sum_{\tau=0}^{D-1}H'(A')^{\tau}\hat{Q}A^{\tau}H
%+\sum_{j=0}^{D}\sigma_j^2C_j' \Big( (A')^{D-j}PA^{D-j} \no\\
%&+\sum_{\tau=0}^{D-1}(A')^{\tau-j}\hat{Q}A^{\tau-j} \Big)C_j \no\\
%&+\sum_{\tau=0}^{D-1} \Big( H'(A')^{\tau+1}PH \Lambda^{-1} H'PA^{\tau+1}H \no\\
%&+\sum_{j=0}^\tau \sigma_j^2C_j(A')^{\tau-j+1}PH\Lambda^{-1}H'PA^{\tau-j+1}C_j
% \Big) \bigg]u_t \no\\
%&+\sum_{\tau=0}^{D-1} \hat{\eta}_{t|t+\tau}' (A')^{\tau+2}PH\Lambda^{-1}H'PA^{\tau+2}\hat{\eta}_{t|t+\tau} \no\\
%&-\sum_{\tau=0}^{D-1} \hat{\eta}_{t|t+\tau-1}' (A')^{\tau+2}PH\Lambda^{-1}H'PA^{\tau+2} \hat{\eta}_{t|t+\tau-1} \no\\
%&+\sum_{\tau=0}^{D-1} 2\hat{\eta}_{t|t+\tau}' (A')^{\tau+2}PH\Lambda^{-1}H'PA^{\tau+1}Hu_{t} \no\\
%&-\sum_{\tau=0}^{D-1}2\eta_t A\Big( (A')^{D}PA^{D}+\sum_{\tau=0}^{D-1} (A')^{\tau}\hat{Q}A^{\tau} \Big) H u_t
%\Bigg\}.
%\eeee
%Consider the part of ${\bf E}(\eta_t'\Omega_1\eta_t)$
Utilizing the fact that $A'PH\Lambda^{-1}H'PA=-P+A'PA+\hat{Q}$ with (\ref{system02}) yields that
\bee
&V_{t}(\eta)-V_{t+1}(\eta) \no\\
=& {\bf E} \Bigg\{  \eta_t' \bigg[ (A')^DPA^D+\sum_{\tau=0}^{D-1} (A')^\tau \hat{Q}A^\tau-(A')^{D+1}PA^{D+1}
\no\\
&- \sum_{\tau=0}^{D-1} (A')^{\tau+1} \hat{Q}A^{\tau+1}
+(A')^D( -P+A'PA+\hat{Q} ) A^D  \bigg]   \eta_t
\no\\
&-u_t'  \bigg[ H'(A')^DPA^DH +\sum_{\tau=0}^{D-1}H'(A')^{\tau}\hat{Q}A^{\tau}H \no\\
&+\sum_{\tau=0}^{D-1}  H'(A')^{\tau}(-P+A'PA+\hat{Q})A^{\tau}H
\no\\
&
+\sum_{j=0}^{D}\sigma_j^2C_j' \Big( (A')^{D-j}PA^{D-j}+\sum_{\tau=0}^{D-1}(A')^{\tau-j}\hat{Q}A^{\tau-j} \Big)C_j \no\\
&+\sum_{\tau=0}^{D-1}\sum_{j=0}^\tau \sigma_j^2C_j(A')^{\tau-j}(-P+A'PA+\hat{Q})A^{\tau-j}C_j\bigg]u_t
\no\\
&- \hat{\eta}_{t|t-1}' A'PH\Lambda^{-1}H'PA \hat{\eta}_{t|t-1} \Bigg\} \no\\
&+ {\bf E} \Bigg\{ {\bf E} \bigg[
\sum_{\tau=0}^{D-1} 2\hat{\eta}_{t|t+\tau}' (A')^{\tau+1} (-P+A'PA+\hat{Q}) A^{\tau}Hu_{t} \no\\
&-\sum_{\tau=0}^{D-1}2\eta_t A\Big( (A')^{D}PA^{D}+\sum_{\tau=0}^{D-1} (A')^{\tau}\hat{Q}A^{\tau} \Big) H u_t| \mathcal{F}_{t-1}\bigg] \Bigg\} \no\\
=&{\bf E} \Bigg\{\eta_t'\hat{Q}\eta_t-u_t' \bigg[H'PH+\sum_{\tau=0}^D \sigma_\tau^2 C_\tau'PC_\tau
\no\\
& +\sum_{\tau=1}^D \sigma_\tau^2 C_\tau\left(\sum_{j=-\tau}^{-1} (A')^{j}\hat{Q}A^{j} \right)C_\tau \bigg] u_t \no\\
&-\hat{\eta}_{t|t-1}' A'PH\Lambda^{-1}H'PA\hat{\eta}_{t|t-1}-2\hat{\eta}_{t|t-1}' A'PHu_{t} \Bigg\} \no\\
%=&{\bf E} \bigg\{\eta_t'\hat{Q}\eta_t+u_t'R u_t
%-u_t'\Lambda u_t -2\hat{\eta}_{t|t-1}' A'PHu_{t} \no\\
%&
%-\hat{\eta}_{t|t-1}' A'PH\Lambda^{-1}H'PA\hat{\eta}_{t|t-1} \bigg\} \no\\
=&{\bf E} \bigg\{\eta_t'\hat{Q}\eta_t+u_t'R u_t-\Big( u_t+\Lambda^{-1}H'PA \hat{\eta}_{t|t-1} \Big)'\Lambda \no\\
&\times
\Big( u_t+\Lambda^{-1}H'PA \hat{\eta}_{t|t-1} \Big) \bigg\}.
\eee
Applying
$u_t^*=-\Lambda^{-1}H'PA \hat{\eta}_{t|t-1}$
leads to
\be
V_{t}(\eta)-V_{t+1}(\eta) = {\bf E} \Big[\eta_t'\hat{Q}\eta_t+u_t'R u_t\Big] > 0.
\ee
By means of Lyapunov stability theory and Lemma \ref{lemma01}, we obtain that system \eqref{system03} is mean-square stabilizable while the stabilizing control input is
\beee
&u_t^*=-\Psi^{-1}L'ZA^{D+1} \hat{\eta}_{t|t-1} \\
=&-\Psi^{-1}L'ZA^{D+1} x_t-\sum_{\tau=1}^D \sum_{j=\tau}^D \Psi^{-1}L'ZA^{D+\tau-j} B_j u^*_{t-\tau}.
\eeee

{\em Necessity.}
Suppose system (\ref{system01}) is stabilizable. When $\hat{Q}>0$, the following deterministic system $(A;\hat{Q}^{\frac{1}{2}})$
\bee
x_{t+1}=Ax_t,~y_t=\hat{Q}^{\frac{1}{2}}x_t\label{observe-system}
\eee
is evidently observable. It follows from Theorem 4.2 in \cite{lilin2016}, there exist unique solutions $P_i,~i\in\mathbb{N}_{[1,D+1]}$ satisfying the following coupled algebraic Riccati equations
\bse\label{couple-riccati}
\bee
P_1=& A'P_1A+A'P_{D+1}A+\hat{Q}, \label{couple-riccati-1} \\
P_i=&-(A')^{i-2} M'\Upsilon M A^{i-2},~i\in \mathbb{N}_{[2,D+1]} \label{couple-riccati-i}\\
M=&\sum_{i=1}^{D+1}H'P_iA, \label{couple-riccati-m} \\
\Upsilon =& R+\sum_{i=1}^{D+1}H'P_iH+\sum_{j=0}^{D}\sigma^2_j H_j' P_1 H_j+\sum_{i=2}^{D+1}\sum_{j=0}^{i-2} \sigma^2_j H_j' P_i H_j  \label{couple-riccati-U}\\
H_j=&A^{-j}C_j, ~j\in \mathbb{N}_{[0,D]}
\eee\ese
where $\sum_{i=1}^{D+1}P_i>0$.
In what follows, we simplify the coupled algebraic Riccati equations in \eqref{couple-riccati} to DDARE \eqref{DARE}.
Let us define $P=\sum_{i=1}^{D+1}P_i>0,$
which implies that $M=H'PA$.
From \eqref{couple-riccati-i}, we have
\be
P_2=-M'\Upsilon M,~P_i=A'P_{i-1} A,~i\in \mathbb{N}_{[3,D+1]}. \label{couple-riccati-i+1}
\ee
Taking the sum from $P_1$ to $P_{D+1}$ on both sides of \eqref{couple-riccati-1} and \eqref{couple-riccati-i+1} yields
\bee
P=&A'(P_1+P_2+\cdots+P_{D+1})A+\hat{Q}+P_2  \no\\
=&A'PA+\hat{Q}-M'\Upsilon M.\label{reduce-riccati-P}
\eee
In this case, we have
\beee
P_1=&P-\sum_{i=2}^{D+1}(A')^{i-2} (P-A'PA-\hat{Q}) A^{i-2} \\
=&(A')^DPA^D+\sum_{i=0}^{D-1}(A')^{i}\hat{Q} A^{i},\\
P_i=&(A')^{i-2} (P-A'PA-\hat{Q}) A^{i-2},~i\in \mathbb{N}_{[2,D+1]},
\eeee
%and
%\bee
%P_1=&P-\sum_{i=2}^{D+1}(A')^{i-2} (P-A'PA-\hat{Q}) A^{i-2} \no\\
%=&(A')^DPA^D+\sum_{i=0}^{D-1}(A')^{i}\hat{Q} A^{i}.
%\eee
The parameter $\Upsilon$ satisfies
\bee
\Upsilon
%=& R+H'PH+\sum_{j=0}^{D}\sigma^2_j C_j' (A')^{D-j}PA^{D-j}C_j\no\\
%&+\sum_{j=0}^{D}\sum_{i=0}^{D-1}\sigma^2_j C_j' (A')^{i-j}\hat{Q} A^{i-j} C_j \no\\
%&+\sum_{i=2}^{D+1}\sum_{j=0}^{i-2} \sigma^2_j C_j'(A')^{i-j-2} (P-A'PA-\hat{Q}) A^{i-j-2} C_j \no\\
=&R+H'PH+\sum_{j=0}^{D}\sigma^2_j C_j' (A')^{D-j}PA^{D-j}C_j\no\\
&+\sum_{j=0}^{D-1}\sum_{i=0}^{D-1-j} \sigma^2_j C_j'(A')^{i} (P-A'PA) A^{i} C_j \no\\
&+\sum_{j=0}^{D}\sum_{i=0}^{D-1}\sigma^2_j C_j' (A')^{i-j}\hat{Q} A^{i-j} C_j \no\\
&-\sum_{j=0}^{D-1}\sum_{i=0}^{D-1-j} \sigma^2_j C_j'(A')^{i}\hat{Q} A^{i} C_j \no\\
=& W_{R,\hat{Q}}+H'PH+\sum_{j=0}^D \sigma_j^2 C_j'PC_j,
%+\sum_{j=1}^D \sum_{i=1}^j\sigma_j^2 C_j (A')^{-i}\hat{Q}A^{-i} C_j,
\eee
which is $\Lambda$ in \eqref{dare-parameter-L}.
This proof is complete.
\ok

\section{Proof of Theorem \ref{delay-margin-theorem}}\label{app03}

\textbf{Proof.}
%The proof is similar to that of Theorem 2 in \cite{tan2015} and we only provide a sketch.
System $[A,C_0;\bar{B},\bar{C}]_{D}$ is stabilizable if and only if there exists a unique positive solution $z_i>0$ satisfying
\be
z_i=\mathcal {R}^{r_i}_{D}(z_i)=a^2_i z_i+{q_i}-a^2_i\bar{b}_i^2 \Psi^{-1}_i z^2_i,\label{DARE-scalar}
\ee
where $q_i>0,$ $r_i>0$, and
%the parameter $\Psi_i$ satisfies
\bee
\Psi_i =&\le(\bar{b}^2_i+ \sigma_0^2 c^2_i a^{2D}_i +\sigma_D^2 \bar{c}^2_i a_i^{2D} \ri) z_i
\no\\
&
+ \le(r_i+ \sum_{h=1}^D\sigma_D^2 \bar{c}_i^2 a_i^{2(D-h)}q_i\ri)>0.
\eee
%where
%\bee
%\phi_i=&\bar{b}^2_i+ \sigma_0^2 c^2_i a^{2D}_i +\sigma_D^2 \bar{c}^2_i a_i^{2D}>0, \\
%\varphi_i=&r_i+ \sum_{h=1}^D\sigma_D^2 \bar{c}_i^2 a_i^{2(D-h)}q_i\geq r_i>0.
%\eee
Then,
DDARE \eqref{DARE-scalar} can be rewritten as a quadratic form
\be
c_{i2}z_i^2+c_{i1}z_i+c_{i0}=0, \label{quad}
\ee
where the parameters are
\bee
c_{i2}=&(a_i^2-1)c_i-a_i^2\bar{b}_i^2,
c_{i1}=(a_i^2-1)c_{i0} +c_i, \no\\
c_{i0}=&r_i+ \sum_{h=1}^D\sigma_D^2 \bar{c}_i^2 a_i^{2(D-h)}q_i, \no\\
c_i=&\bar{b}^2_i+ \sigma_0^2 c^2_i a^{2D}_i +\sigma_D^2 \bar{c}^2_i a_i^{2D} \no.
\eee
Similar to the proof of Theorem 2 in \cite{tan2015}, the above quadratic equation \eqref{quad} has a unique positive solution if and only if $c_{i2}<0$, i.e.,
%the following algebraic inequality holds,
\be
(a_i^2-1)\le( \bar{b}^2_i+ \sigma_0^2 c^2_i a^{2D}_i +\sigma_D^2 \bar{c}^2_i a_i^{2D} \ri)
-a_i^2\bar{b}_i^2<0,
\ee
which is equivalent to
\be
a^{2D}_i (\sigma_0^2 c^2_i+\sigma_D^2 \bar{c}^2_i)(a_i^2-1)<\bar{b}^2_i.\label{algebraic-inq}
\ee
When $a_i^2\leq 1$, inequality \eqref{algebraic-inq} holds for any bounded delay $D_i>0$, and hence, sub-system $[a_i,c_i;\bar{b}_i,\bar{c}_i]_{D}$ is stabilizable for $D<D^i_{\max}=+\infty$.
When $a_i^2> 1$, it follows that
% from \eqref{algebraic-inq} that
\be
a_i^2-1<a^{2D}_i (a_i^2-1)< \frac{\bar{b}^2_i}{\sigma_0^2 c^2_i+\sigma_D^2 \bar{c}^2_i}=h_i.
\ee
The delay margin of $[a_i,c_i;\bar{b}_i,\bar{c}_i]_{D}$ is derived as
\beee
D^i_{\max}=\left\{
\begin{array}{lr}
\frac{\ln(h_i)-\ln(a_i^2-1)}{2\ln(a_i)},&{\rm if}~1<a_i^2<h_i+1, \\
0,&{\rm if}~a_i^2\geq h_i+1.
\end{array}
\right.
\eeee
%On the other hand,
%the delay margin $D_{\max}$ of system $[A,C_0;\bar{B},\bar{C}]_{D}$ satisfies
%\be
%D_{\max}=\min_{1\leq i\leq n}D^i_{\max}.
%\ee
To sum up, if $a_i^2=1$ holds for each $a_i\in U_A$, system $[A,C_0;\bar{B},\bar{C}]_{D}$ is stabilizable for $D\geq 0$,
and if $a_i^2 \geq h_i+1$ for some $a_i\in U_A$, system $[A,C_0;\bar{B},\bar{C}]_{D}$ cannot be stabilized for any $D> 0$.
Otherwise, we obtain that
\be
D_{\max}=\min_{a^2_i>1}\frac{\ln(h_i)-\ln(a_i^2-1)}{\ln(a_i^2)},~a_i\in U_A,
\ee
which completes the proof.
\ok
%%%%%%%%%%%%%%%%%%%%%%%%%%%%%%%%%%%%%%%%%%%%%%%%%%%%%%%%%%%%%%%%%%%%%%%%%%%%%%%%%%%%%%%%%%%%%%%%%%%%%%%%%%%%%%%%%%%%%%%%%%%%%%%%%%%%%%%%%%%%%%%%%%

\end{document}